# BALANCED SCHEDULING OF SCHOOL BUS TRIPS USING A PERFECT MATCHING HEURISTIC

**Ali Shafahi (Corresponding Author)**
PhD Candidate
Department of Computer Science
University of Maryland
College Park, MD 20742
Email: ashafahi@umd.edu

**Sanaz Aliari**
PhD Candidate
Department of Civil and Environmental Engineering
University of Maryland
College Park, MD 20742
Email: saliyari@umd.edu

**Ali Haghani**
Professor
Department of Civil and Environmental Engineering
University of Maryland
College Park, MD 20742
Email: haghani@umd.edu

June 27th, 2018

## ABSTRACT

In the school bus scheduling problem, the main contributing factor to the cost is the number of buses needed for the operations. However, when subcontracting the pupils' transportation, unbalanced tours can increase the costs significantly as the lengths of some tours can exceed the daily fixed driving goal and will result in over-hour charges. This paper proposes an MIP model and a matching-based heuristic algorithm to solve the "balanced" school bus scheduling problem with fixed start times in a multi-school setting. The heuristic solutions always have the minimum number of buses as it starts with a minimal number of tours and does not alter the number of tours during its balancing stage. The effectiveness of the heuristic is tested by comparing its solutions with results from solving the MIP using commercial solvers whenever solvers could find a good solution. To illustrate the performance of the MIP and the heuristic, 11 problems were examined with different numbers of trips which are all based on two real-world problems: a California case study with 54 trips and the Howard County Public School System with 994 trips. Our numerical results indicate the proposed heuristic algorithm can find reasonable solutions in a significantly shorter time. The balanced solutions of our algorithm can save up to 16% of school bus operation costs compared to the best solution found by solvers from optimizing the MIP model after 40 hours. The balancing stage of the heuristic decreases the Standard Deviation of the tour durations by up to 47%.

## INTRODUCTION

Optimal scheduling of bus services can lead to significant savings in operating costs. An example of such a service exists in the school-bus service of Howard County Public School System (HCPSS) in the state of Maryland, which serves 78 schools with a total of 994 afternoon-trips. Each afternoon trip starts from a school and finishes at the last passenger's delivery point. Planning and implementation of cost efficient school-bus transportation systems are indeed a challenging task for such a large demand network. It is noteworthy that the bell-times of these schools are different and hence, it is possible to form links of sets of trips into blocks/tours that each can be served by a single bus. The problem of chaining trips, each being served by a single bus is called school bus blocking/scheduling problem. Optimal bus blocking in a way that the number of blocks (and hence the number of required buses) is minimized can significantly reduce the fixed operating costs of a school district. Considering the very high annual fixed cost for each bus (between $50,000 and $100,000 in the state of Maryland), the minimization of the total number of buses is a major contributor to cost efficiency of school bus operations. The variable component of cost, which is proportional to the traveled time/distance of buses, is the second largest contributor to costs and should also be minimized. Given that in blocking/scheduling problems, the trip durations are fixed, the latter can be achieved by minimizing the aggregate portions of operation times of each bus, during which they have no passengers (deadhead). Figure 1 visualizes an example of a small school-bus system and one of the corresponding possible blockings.

In this paper, a variation of the scheduling/blocking problem for school buses is introduced that can further improve the cost efficiency of the service. More precisely, the considered problem is the balanced school bus scheduling problem (SBSP) in which the objective is not only to seek minimization of the total number of buses but also to provide comparable (balanced) block/tour durations for all the buses. The balanced blocking is beneficial because more savings can be gained by reducing the over hour operating costs and avoiding extra charges. Moreover, such a balanced blocking would be more equitable for the drivers. This problem is formulated as a Mixed Integer Programming (MIP) model. However, since finding the exact solution of such a model is an NP-hard problem, a two-stage heuristic approach is proposed to find a balanced blocking solution for large-size problems. First, the generic SBSP with fixed trip start times (without the balancing constraints) is solved to find an initial bus blocking schedule that simultaneously minimizes the number of blocks and the aggregate deadhead. Then, the schedules are become more balanced by building a bipartite graph and using a perfect matching approach.

The remaining of the paper is organized as follows. The next section presents a brief literature review of the school bus scheduling problem. Then a detailed description of the problem is provided. Subsequently, the MIP and the matching-based heuristic are introduced. Afterwards, the results from the mathematical model and the heuristic algorithm are presented and compared and the effects of the balancing stage of the heuristic on the set of tours are illustrated. Finally, the conclusions and the possible future extension of the studied problem are presented.

## LITERATURE REVIEW

The school bus blocking/scheduling problem (SBSP) belongs to the general class of Open Vehicle Routing Problems (OVRP) (1). The goal of the OVRP in general is to design a set of Hamiltonian paths (open routes) to serve customer demand. In the case of school bus scheduling, the problem consists of chaining a set of scheduled trips to an optimal number of blocks that are bound by the number of available vehicles. Moreover, there are more characteristics to the SBSP that differentiates it from the generalized version of OVRP such as consideration of heterogeneous fleet, mixed student load, multiple depots; and work-load balancing for the school bus trips. The OVRP has been extensively studied in the literature. A survey on the algorithms for solving OVRP is provided in (2). These studies differ in the characteristics of the problem, the objective function and constraints used to develop the mathematical model and the strategy used to find the solution for the model.

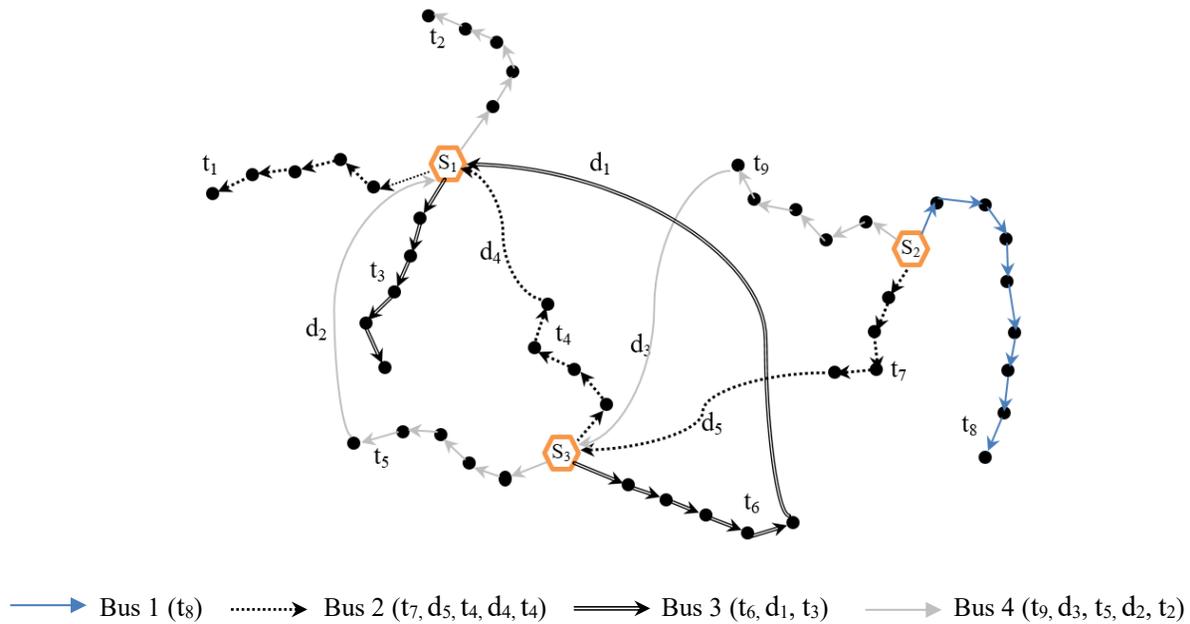

Bus 1 (t₈) Bus 2 (t₇, d₅, t₄, d₄, t₄) Bus 3 (t₆, d₁, t₃) Bus 4 (t₉, d₃, t₅, d₂, t₂)

**Figure 1 Illustration of an unbalanced blocking for an afternoon 9-trip, 3-school scenario using four buses**

As for the characteristics of the SBSP, many different aspects of a practical transportation system are considered. For instance, some of the earliest research consider single depot for buses (3, 4), while some others take the possibility of having multiple depots into account (5 - 7). Some of the more practical aspects of the bus scheduling problem such as fuel consumption (5), fixed/flexible time windows (6, 7) and multiple vehicle types (8, 9) have also been taken into account.

The SBSP might then be approached through different practical objectives; one might be solely interested in maximizing profit (minimizing costs), while others take measures of the equity of service (fairness) and level of service into account. The most common objective functions are the minimization of total travel distance or time, the total number of trips, and the total cost that is the combination of the aggregate travel distance and number of trips (10-14). Setting aside the fixed costs for each trip, the bus operational/variable cost is directly related to the total travel distance. Multi-objective cost functions, composed of conflicting objectives, are also widely used. In (14-15), the two conflicting objectives of minimizing the number of buses while minimizing the maximum ride time for riders are considered. These are conflicting objectives since fewer buses would require the routes to be longer. A similar objective function is used in (16) where the authors formulate the blocking problem as an OVRP with fixed time-windows (by considering each trip as a virtual stop) and present two exact solution algorithms based on the Assignment Problem (AP). One for the special case where the start and end service time of each virtual stop are fixed, and the buses are homogenous, and another AP based solution algorithm for the more general case. The second approach is based on an iterative heuristic algorithm, composed of an AP based construction algorithm and an improvement step.

Balancing of school transportation systems has been studied in the past. However, these studies have all focused on balancing the routes/trips and not the tours which are a collection of the trips. Balanced route/trip length is sometimes used in multi-objective cost functions and is motivated by the need for equity and fairness of the service (10). In (10), some equity criteria including balanced route/trip length for different buses are considered in the form of a multi-objective nonlinear integer program. The authors present a heuristic method based on an initial clustering of students into fixed clusters of routes. Although the heuristic did not explicitly enforce the balancing restrictions, the authors found that the durations of the routes are close to each other. Load balancing is another example of balancing in the school transportation operations. The goal of load balancing is to have routes which have a similar number of students. A case study of the load balancing of school bus routing in Hong Kong is presented in (17). They use a combinatorial multi-objective optimization approach and addresses the travel length by a modification stage to the initial solution, during which the spare capacity of buses with empty spaces are filled with some of the load of fuller buses that serve similar routes.

Matching-like algorithms have been previously used in the school bus transportation context. Perfect matching is used in (18) to assign a morning route to an afternoon route and construct work shifts for bus drivers. In (19) a repeated matching heuristic is introduced to solve the bus routing problem, to find optimal clustering of customers that can be served by each bus. For balancing purposes, matching has also been used. However, to the best of the authors' knowledge, it has not been used for balancing of school bus tours. In (20), a relaxation of the matching problem, called semi-matching, is described and applied in the context of a job assignment system, aiming at balancing the load assigned to each machine. This approach cannot be directly applied to our problem since semi-matching cannot assure a compatible blocking when deadheads are considered.

## PROBLEM STATEMENT

This study presents a mathematical model and a solution heuristic for the balanced school bus scheduling (tour generation) task for a school district with multiple schools. The inputs to this problem are the set of afternoon trips (starting from school) and their information, and a goal hour which is used for balancing and also as a soft cut-off for the maximum duration/mileage of any single tour (block). The input trips (their sequence of stops, duration, and start-time) are fixed and will not be altered. The goal of the balanced scheduling problem is to assign the trips to tours such that primarily the number of tours are minimized, and secondarily the tours are mostly shorter than a given goal and hopefully balanced. The main objective is minimizing the number of tours since as previously mentioned the number of tours is the main contributor to the transportation costs of school districts as each additional tour requires the deployment of an additional bus for a school year.

The secondary objective is to ensure that the aggregate duration that tours exceed the given goal is minimized. This objective is particularly useful for the cases that the transportation services of a school district are outsourced. Usually, each contracted bus has a fixed mileage/time for its service and any additional service beyond that is charged extra at a high over-time rate. In this study, the tour duration defined to be the student ride time (trip duration) plus the deadhead (a portion of operational times with no students on board).

It is worth noting that fixing the school trips are justified since many of the large school districts such as the Philadelphia school district continue to manually build the school trips themselves via help from their schedulers because of many local and social constraints that are needed to be taken into account when building every single trip. Examples of these constraints are parents' specific requests, safety considerations based on the crime level of the school district, possibilities for students to cross streets for different types of streets, etc.

An assumption implied throughout this study is that all trips can be assigned to all tours (homogenous fleet). This assumption can be relaxed by solving different problems for a different set of trips all having similar loads. Although this approach may not yield the optimal solution, it can lead to find an approximate solution. However, it is worth noting that the homogeneous fleet assumption is not very limiting as it is valid for many school districts which outsource their logistics. Our real-world instance case, HCPSS, is an example which outsources all of their regular education buses and uses homogeneous fleet.

## METHODOLOGY

In this section, a mathematical model that is used for finding the optimal solution of the balanced school bus blocking problem is presented. Then, a solution heuristic algorithm that can be used for solving large instances of the problem is introduced.

### Mixed integer programming model

The set $T = \{1, 2, \ldots, N\}$ of trips are given to serve the demand within the school district. It is desirable to link the trips into the minimum K tours to hopefully reduce the number of vehicles needed for operations. Since each tour is operated by one bus, the set $V = \{1, 2, \ldots, K\}$ denotes both the set of vehicles and tours which have more than one trip assigned to them. Note that since there are N trips, K would be at most N/2. Consequently, $K = N/2$ is set. This way of defining tours reduces the number of variables which are related to both trips and tours since the number of potential tours are cut into half by only considering tours that have more than one trip assigned to them. Due to the homogeneous-fleet assignment, the trips which are done just by themselves would not need to be assigned to a particular tour. The other inputs to the problem include the trips' travel times, the fixed start time and the pairwise deadhead matrix between all pairs of the compatible trips. A pair of trips are compatible if it is possible to chain them in a tour without violating the time window constraint. The following notations are defined to model the balanced school bus blocking problem.

| |
|---|
| $x_{ij}^k$: 1 if trips $i$ and $j$ are assigned to tour $k$ and trip $i$ is served right before trip $j$; 0, otherwise. |
| $m_i^k$: 1 if trip $i$ has at least one trip preceding and at least one trip succeeding it in tour $k$; 0, otherwise. |
| $n_i^k$: 1 if trip $i$ is served by tour $k$; 0, otherwise. |
| $a_i$: 1 if trip $i$ belongs to a single-trip tour (no trips succeeding or preceding it) |
| $b^k$: 1 if tour $k$ is being used to serve at least two trips; 0, otherwise |
| $p^k$: The period of time that tour $k$ exceeds the goal |

The notations for the inputs of our problem are listed as:

| | | |
|---|---|---|
| $T_i$ | : | Travel time duration for trip $i$. This is equal to the student ride time for trip $i$. |
| $D_{ij}$ | : | Deadhead travel time between trips $i$ and $j$. |
| $G$ | : | Maximum goal duration |
| B | : | Maximum number of trips that can be assigned to a single tour |
| $M_{bus}$ | : | Penalty term for dispatching each additional bus |
| $M_G$ | : | Penalty term for exceeding the set goal duration $G$. |

The objective function (1) aims at minimizing the total penalty costs associated with dispatching additional vehicles to satisfy the demand ($M_{bus}$) and exceeding the goal ($M_G$) for a scheduling plan.

$$Minimize\ z = M_{bus} \cdot Nbuses + \sum_k M_G \cdot p^k \qquad (1)$$

$$Nbuses \coloneqq \sum_k \sum_i \sum_j x_{ij}^k - \sum_k \sum_i m_i^k + \sum_i a_i$$

This objective function is subject to a number of constraints. To find the penalty associated with the overtime operations, the duration of time that exceeds the set goal must be obtained. Constraints (2) calculate the amount of time that each generated tour exceeds the set goal. Which means a value is assigned to $p^k$ which acts as a slack variable during the calculations.

$$\sum_i \sum_j D_{ij} x_{ij}^k + \sum_i \sum_j (T_i + T_j) x_{ij}^k - \sum_i T_i m_i^k - p^k \le G \qquad \forall k \in V \qquad (2)$$

Constraints (3) enforce that all trips must be assigned to exactly one tour. Either a tour that is only composed of one trip (when $a_i = 1$) or a tour which has more than one trip assigned to it.

$$\sum_j \sum_k x_{ij}^k + \sum_j \sum_k x_{ji}^k + a_i - \sum_k m_i^k = 1 \qquad \forall i \in T \qquad (3)$$

Constraints (4) and (5) find non-zero $m_i^k$ which indicates the middle trips.

$$\sum_j x_{ij}^k + \sum_j x_{ji}^k \ge 1 + m_i^k - (1 - n_i^k) \qquad \forall i \in T, \forall k \in V \qquad (4)$$

$$\sum_j x_{ij}^k + \sum_j x_{ij}^k \le 2 \times n_i^k \qquad \forall i \in T, \forall k \in V \qquad (5)$$

Each generated tour must have exactly one starting trip and exactly one ending trip that are not middle trips.

$$\sum_i n_i^k - \sum_i m_i^k = 2 \times b^k \qquad \forall k \in V \qquad (6)$$

Constraints (7) find non-zero $a_i$ which indicates that trip $i$ belongs to a tour with a single trip.

$$\sum_j \sum_k x_{ij}^k + \sum_j \sum_k x_{ji}^k \le 2(1 - a_i) \qquad \forall i \in T \qquad (7)$$

Constraints (8) and (9) find non-zero $b^k$ which indicates the non-single tours:

$$\sum_i \sum_j x_{ij}^k \leq (B-1) \cdot b^k \qquad\qquad \forall k \in V \qquad (8)$$

$$\sum_i \sum_j x_{ij}^k \geq b^k \qquad\qquad \forall k \in V \qquad (9)$$

Each trip can appear at most once as a preceding trip and once as a succeeding trip:

$$\sum_j \sum_k x_{ij}^k \leq 1 \qquad\qquad \forall i \in T \qquad (10)$$

$$\sum_j \sum_k x_{ji}^k \leq 1 \qquad\qquad \forall i \in T \qquad (11)$$

Constraints (12) reduce the symmetry of the problem resulted by homogenous fleet:

$$b^k \leq b^{k-1} \qquad\qquad \forall k \in \{2, \ldots, N\} \qquad (12)$$

And finally binary and non-negativity constraints:

$$a_i, m_i^k, n_i^k, b^k, x_{ij}^k \in \{0,1\} \qquad\qquad \forall i, j \in T, \forall k \in V \qquad (13)$$

$$p^k \geq 0 \qquad\qquad \forall k \in V \qquad (14)$$

This formulation is based on the work in (**Error! Reference source not found.**) with some modifications to address the balancing characteristics. It is important to note that since this formulation is based on preprocessing the data, the time window and sub-tour elimination constraints are not needed. Indeed, time window constraints are implicitly considered by excluding incompatible trips from the input data. For the same reason, not many sub-tours are possible: for instance, if trip i is compatible with trip j, trip j cannot be compatible with trip i since one of them starts later than the other.

**Perfect matching heuristic approach:**

As mentioned before Open Vehicle Routing Problems (OVRP) are categorized as NP-hard problems, and it is impossible to obtain optimal solutions for large networks. Indeed, in our case, the commercial solver used (Xpress) for solving the proposed mathematical model yields a suboptimal solution with a 15% optimality gap after 42 hours for a problem with 250 trips. Therefore, a heuristic approach to tackle the large size problems is proposed. The main idea of this heuristic is based on breaking the problem into two stages: first, in the ***blocking*** stage, an optimal solution will be obtained for the blocking problem using the assignment-based approach discussed in (16) which applies to the blocking problems with fixed start times. The first stage generates solutions that have the minimum number of tours, and the aggregate deadhead is also minimized. Second, in the ***balancing*** stage, the solution will be modified to satisfy the duration goal as much as possible. The modification is done through finding a minimum weighted perfect matching in a bipartite graph. The Minimum weight perfect matching problem is also called the assignment problem and can be solved in polynomial time using different algorithms such as Hungarian algorithm, modified primal simplex, and the auction algorithm. Since the balancing step is based on a perfect matching of tours to the trips, the number of tours that are minimal from the blocking step will not be altered. The perfect matching step ensures that the heuristic always gives the minimal number of tours and the main fixed-cost contributor (Number of buses) is minimal. The steps of this heuristic are described below, and a summary of the steps is depicted in Figure 2.

- **Stage 1:** *Optimal blocking:* Relax the time balancing constraint and obtain an optimal schedule for the input trips using the assignment-based approach described in (16).
- **Stage 2:** *Balancing modification via perfect matching*:
    - **Step 2-1.** *Bipartite graph vertex construction:* Build the bipartite graph as follows: Strip the last trip from each tour, and associate a vertex with each stripped tour, and build set A.

      Similarly, associate a vertex with each stripped trip and build set B. Two sets of vertices of equal size are obtained ready for the bipartite graph.
    - **Step 2-2.** *Compatible edge construction:* Add an edge between set A, the stripped tours, and set B, stripped trips, in the bipartite graph, if the linked pair has a compatible* time window. Afterward, calculate the corresponding weights**.

      \* Trips a and b are compatible if the end time of current trip (trip a) plus the deadhead to the next trip ($D_{ab}$) will not exceed its (trip b) start time.

      \** The weight for a newly constructed tour is the amount of time the tour operates beyond the goal ($p^{knew}$).
    - **Step 2-3.** *Minimum weight perfect matching:* Perform the Hungarian algorithm (Kuhn–Munkres algorithm) to find the optimal matching.

Granted that the problem is broken into two assignment problems, a drastic decrease in the computation time is expected.

## COMPUTATIONAL RESULTS

In this section, the results of both the mathematical model and the heuristic algorithm are discussed for 11 different scenarios. All scenarios are built based on two real-world problems. Nine scenarios are built by sampling from the trips of a real-world case study. For these nine scenarios, a scenario with $n$ trips is built by taking $n$ ($n$ =10,20,30,40,50,100,200,250,300,500) random trips out of the maximum 994 trips of the real-world data of Howard County Public School System (HCPSS). The Scenario names are set such that they reflect the number of trips associated with them. The average duration of the trips is around 25 minutes with the standard deviation of 10 minutes. The minimum and the maximum duration of the HCPSS trips are 7 and 84 minutes, respectively. A smaller-size real world problem is also considered for a school district in California that has 54 trips (averages duration: 40 minutes, minimum duration: 22 minutes, maximum duration: 66 minutes, and the SD: 10 minutes). Given the distribution of the duration for trips and deadheads, it was noticed that a 75-minute goal would illustrate the balancing effects best thus, the goal is set to be equal to 75 minutes in all problems.

The MIP problem is solved by Xpress on a dual processor Intel Xeon computer with a total 192 GB of memory and 40 cores each 2.20GHz. Xpress uses a combination of branch and bound and heuristic searches at the nodes to search for the optimal solution. Xpress reported the optimality of the solution for all scenarios with up to 50 trips. However, for larger problems, it either was not able to load the problem on memory or was not able to report optimality even after 40 hours.

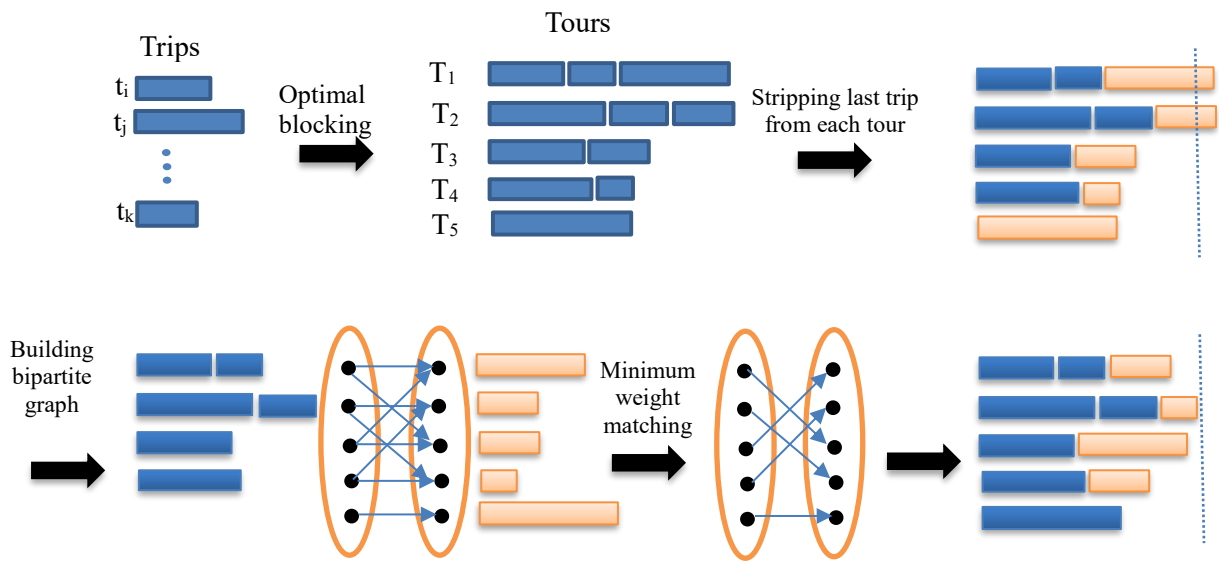

**Figure 2 Schematic representation of the proposed matching-based heuristic**

To see the effect of the coefficients of the objective function on the performance, the model is run using two values for the penalty of additional buses, $M_{bus}$. A larger $M_{bus}$ (60,000) and a smaller $M_{bus}$ (600). The larger value is the average cost for each additional bus in Maryland. The smaller value is set such that the cost of adding a new bus would never exceed the total cost of exceeding the goal. Figure 3 illustrates various performance indicators for mainly the scenarios in which Xpress was able to find a solution. As it can be seen in the results (Figure 3), the performance of the solutions found by solving the MIP using Xpress is sensitive to the coefficient of buses. This behavior is originated from the rounding errors and the precision of the machine. When the lower bound of the objective function is a large number, say in the order of 1 million, then an objective function value of 1,000,001 has a gap of %1e-6 and is practically optimal and is less than the default tolerance of Xpress (%1e-4).

The maximum running time for the scenarios with less than 200 trips is set to 10 hours and the scenarios exceeding 250 trips to 40 hours. For 300 trips, no valuable solution was found even after 40 hours as the problem size is so big that loading the problem is almost impossible. The number of variables, constraints, and nonzero elements for the 300-trip scenario is roughly 3.68 million, 77.21 thousand, and 29.87 million respectively. Each simplex iteration takes 11,743 seconds. Letting Xpress run for 150,000 seconds will only give a solution that requires 300 buses (each trip assigned to its exclusive tour). The problems with 500 trips and 994 trips were not loaded on the computer after 150,000 seconds.

**Evaluation of the heuristic by comparing with good solution of the MIP:**

To compare the solution of the proposed matching heuristic, the outputs of the heuristic are compared with the solutions found by solving the MIP using Xpress. Xpress uses a combination of branch and cut and local search heuristics for solving an MIP. For smaller sized problems,

Xpress can find the optimal solution. One important measure and contributor to the overall transportation cost is the number of buses/tours. The number of buses for the different scenarios are depicted in Figure 3-a. We can see that for cases in which the problem size is small/medium (i.e. the number of trips are less than 100), the Xpress solver can find the minimum number of buses. However, as the problem size exceeds 100 trips, the MIP solver is no longer able to find the minimum number of buses during its set maximum run-time. For the scenarios with 500 trips or more in which the MIP solver was not loaded within a 150,000 second time limit, the number of tours was set equal to the number of trips. Based on these results, it is observable that the MIP solver loses its value very soon and is only practical for small to medium size scenarios. Thus, it cannot be used for the HCPSS problem as it will not give any "good" solutions within any foreseeable time. However, it can find the least number of buses for the small size California problem.

The second contributor to cost is the excess minutes (Figure 3-b). The transportation authorities have to pay extra hour fees for those. Comparing the aggregate over time minutes of the heuristic solutions with the MIP solutions for all scenarios is not fair since in many of the scenarios the Xpress was not able to find a solution with minimum number of tours within the set time limit. At the extreme case, if each bus just makes one trip then there will be no excess minutes, but obviously, the overall cost is a lot more. To make the comparison fair, it is only focused on the scenarios in which the MIP solver (Xpress) was able to find the minimum number of tours (N10, N20, N30, N40, N50, N100, N54) or its number of tours was very close to the minimum number (N200). As it can be seen in Figure3-b, the heuristics solutions are comparable with the solutions from the MIP solution. A better measure to compare the solutions of the heuristic with those of the MIP might be the annual cost of each solution that combines the number of buses and the over-hour costs.

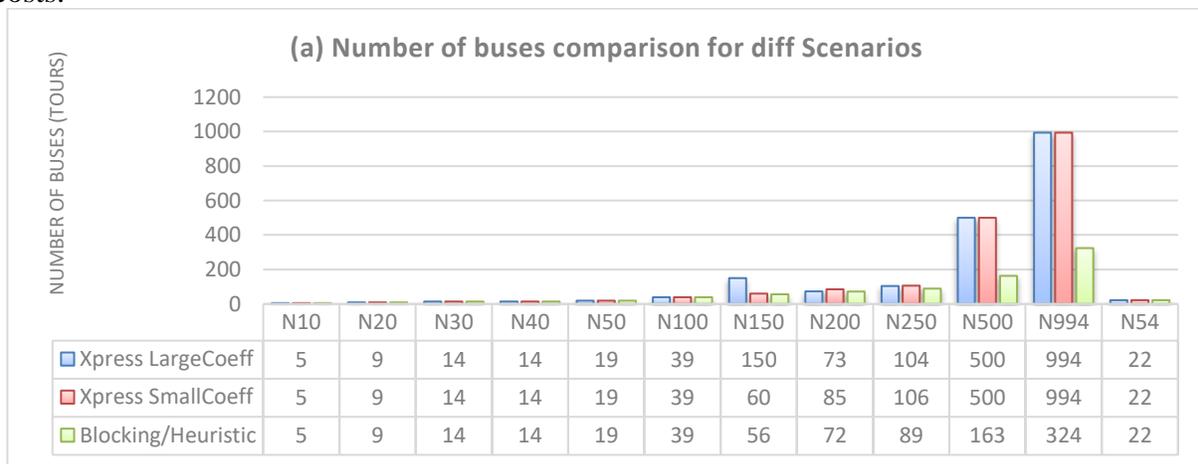

**(a) Number of buses comparison for diff Scenarios**

| | N10 | N20 | N30 | N40 | N50 | N100 | N150 | N200 | N250 | N500 | N994 | N54 |
|---|---|---|---|---|---|---|---|---|---|---|---|---|
| Xpress LargeCoeff | 5 | 9 | 14 | 14 | 19 | 39 | 150 | 73 | 104 | 500 | 994 | 22 |
| Xpress SmallCoeff | 5 | 9 | 14 | 14 | 19 | 39 | 60 | 85 | 106 | 500 | 994 | 22 |
| Blocking/Heuristic | 5 | 9 | 14 | 14 | 19 | 39 | 56 | 72 | 89 | 163 | 324 | 22 |

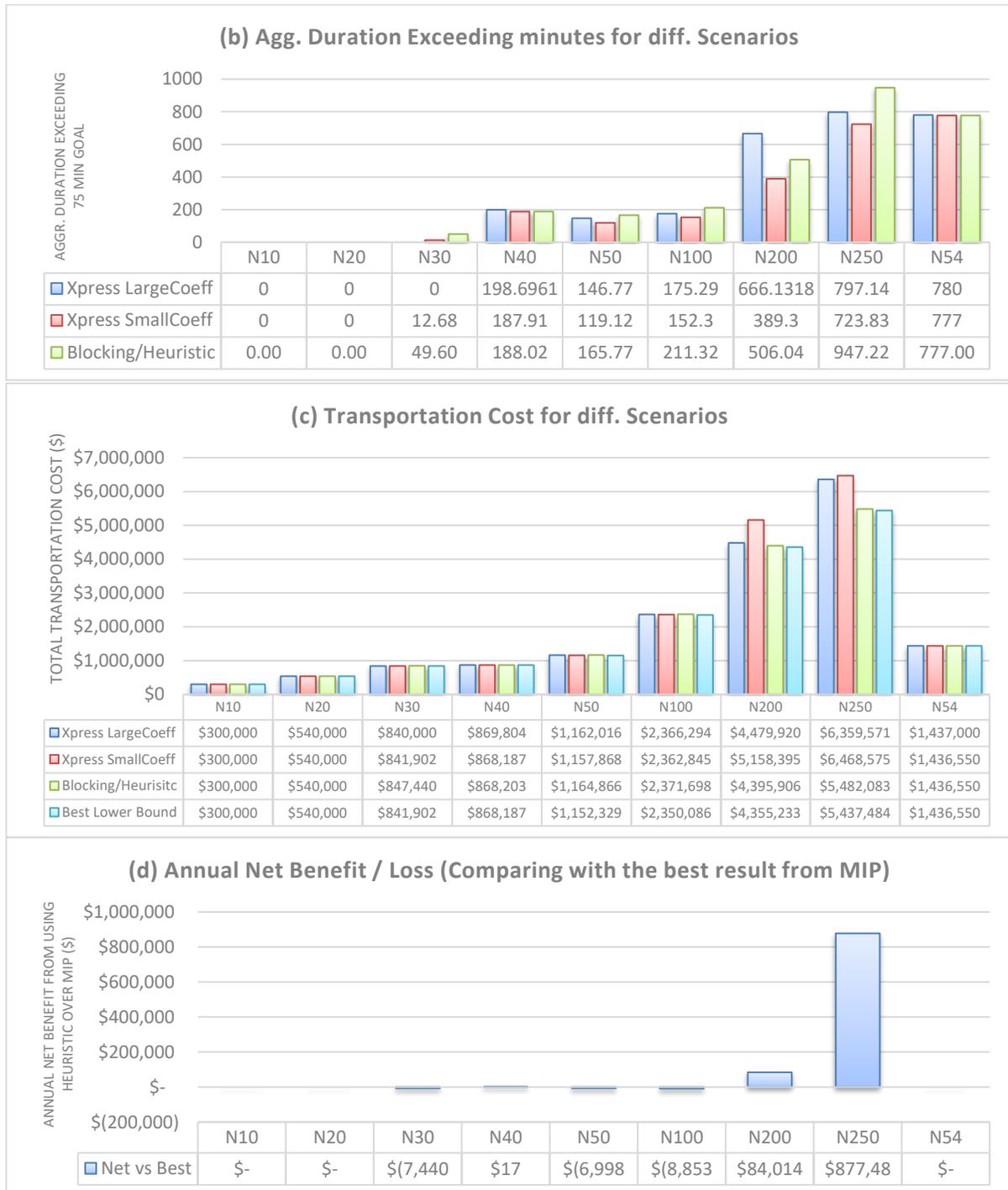

**Figure 3 (a) Number of buses needed; (b) aggregate minutes exceeding the 75-minute goal; (c) annual transportation costs; and (d) The net annual benefit/cost of using the heuristic compared with the best solution found by solving the MIP using Xpress for different scenarios.**

Assuming an extra-hour rate of $50/hr and 180 school days per year, the cost for each additional minute per year is $150. Also, assuming $60,000, the overall cost of the solutions gained from the heuristic and MIPs as well as the best lower bound found by Xpress are summarized in

Figure 3-c. Comparison of the best lower bounds of the MIP found by Xpress and the heuristic solutions indicates a maximum gap of 1.1% and an average gap of 0.5% which accentuates the effectiveness of the proposed heuristic. The net annual benefit/cost of using the heuristic over the best result from the MIP is depicted in Figure 3-d. As it can be seen, for those scenarios which are small enough that the MIP solver can find a good solution, the additional annual cost of the heuristic's solution is always less than $10,000. In the California example (N54), the heuristic finds the optimal solution. However, as the problem size increases, the savings from the heuristic increase. The heuristic can find a solution that has about 16% less cost for the scenario with 250 trips. For larger scenarios, the MIP solver cannot find any reasonable solution at all even when it is given a 150,000 second time limit.

Another measure that shows how well the heuristic performs is the number of tours exceeding the goal. This measurement is interesting as over-hour costs have to be paid for each of the tours exceeding the goal. If there is a fixed cost component to over hour-costs (such as having a one-hour minimum duration for over-hours), this measurement turns into a significant one. For example, if two trips exceed the goal by 5 minutes, the transportation district might need to purchase an additional hour for each of the buses due to contractual obligations. While this measurement is not optimized, it is still exciting to see that the heuristic performs okay on this measurement. As illustrated in Figure 4, the heuristic' performance is close to the MIP solutions in all scenarios except for scenario N250 for a 75 minute goal.

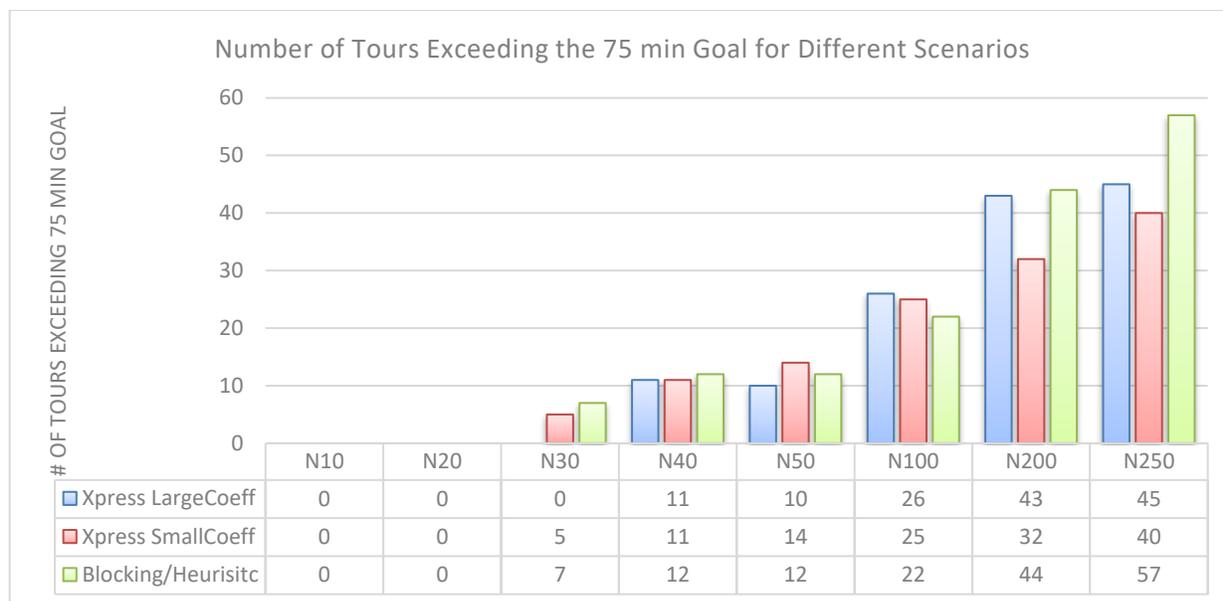

**Figure 4 Comparison of different methods regarding number of tours exceeding the goal**

**How the heuristic affects the schedule's key performance indicators (KPI):**

The number of buses is a critical KPI. Due to the optimality of the solution during the blocking stage, it is known that the heuristic provides the minimum number of buses. To illustrate the effect of the balancing step, the optimal solutions of the blocking problem are compared with the new solution from the matching heuristic. Recall that the optimal solutions to the blocking problem primarily minimize the number of tours (buses) and secondarily minimize the aggregate deadhead time. The heuristic's performance can be measured based on statistics regarding the tour durations before and after applying the second stage (balancing-stage). In particular, the focus is on three KPIs: a) standard deviation of the tours (SD); b) aggregate minutes exceeding the set goal, and c) percentage of tours exceeding the set goal before (old) and after (new) the balancing stage. Also, three statistics about the tours are reported, namely, their average, minimum tour duration, and maximum tour duration. These results are summarized in Table 1 and Table 2. Table 1 shows that the heuristic acts by increasing the minimum tour duration and usually decreasing the maximum tour duration. This result is making the tours more balanced as it can be seen in Table 2. The only case in which the standard deviation increases after applying the heuristic is the California case study (N54). This increase is only a marginal increase. Therefore, it is observable that the heuristic's performance is very satisfactory regarding making more balanced tours. The first two scenarios with 10 and 20 trips remain unchanged after the balancing step because the duration of their tours is not over the 75-minute goal and hence all of the weights of the matching are equal to zero, and the balancing step does nothing.

**Table 1 How the heuristic changes the tours**

| Scenario | Perfect matching heuristic | | | | | | Heuristic running time(s) |
|---|---|---|---|---|---|---|---|
| | **Statistics** | | | | | | |
| | Min(old) | Min(new) | Max(old) | Max(new) | Avg(old) | Avg(new) | |
| N10 | 37.32 | 37.32 | 71.00 | 71.00 | 53.26 | 53.26 | 0 |
| N20 | 45.72 | 45.72 | 72.14 | 72.14 | 60.04 | 60.04 | < 0.01 s |
| N30 | 36.46 | 46.98 | 109.63 | 91.86 | 68.89 | 72.47 | < 0.01 s |
| N40 | 62.15 | 69.70 | 112.54 | 112.54 | 87.76 | 87.91 | <0.01 s |
| N50 | 41.10 | 47.34 | 112.23 | 112.23 | 78.37 | 80.73 | 0.01 |
| N100 | 42.31 | 51.72 | 124.02 | 113.74 | 75.58 | 78.54 | 0.12 |
| N200 | 33.15 | 44.82 | 110.89 | 120.89 | 77.59 | 79.70 | 0.5 |
| N250 | 37.78 | 51.19 | 170.02 | 128.04 | 82.48 | 83.49 | 0.64 |
| N500 | 22.64 | 70.77 | 145.58 | 150.31 | 86.89 | 87.20 | 2.94 |
| N994 | 47.70 | 61.56 | 153.68 | 148.73 | 87.87 | 88.20 | 13.99 |
| N54 | 69.00 | 76.00 | 145.00 | 145.00 | 110.32 | 110.32 | 0.01 |

**Table 2 Key Performance Indicators of the scenarios before and after applying the heuristic**

| Scenario | Perfect matching heuristic | | | | | |
|---|---|---|---|---|---|---|
| | **KPI** | | | | | |
| | SD(old) | SD(new) | Rat. Exceed. tours(old) | Rat. Exceed. tours(new) | Exceeding mins(old) | Exceeding mins(new) |
| N10 | 13.18 | **13.18** | 0/5 | 0/5 | 0.00 | 0.00 |
| N20 | 9.11 | **9.11** | 0/9 | 0/9 | 0.00 | 0.00 |
| N30 | 20.83 | **12.25** | 6/14 | 7/14 | 88.86 | 49.60 |
| N40 | 14.44 | **13.15** | 11/14 | 12/14 | 198.70 | 188.02 |
| N50 | 18.97 | **15.16** | 12/19 | 12/19 | 182.59 | 165.77 |
| N100 | 19.71 | **10.45** | 19/39 | 22/39 | 321.08 | 211.32 |
| N200 | 18.79 | **12.52** | 43/72 | 44/72 | 660.66 | 506.04 |
| N250 | 21.17 | **15.69** | 51/89 | 57/89 | 1077.75 | 947.22 |
| N500 | 18.08 | **13.40** | 126/163 | 184/163 | 2319.05 | 2008.42 |
| N994 | 16.94 | **13.43** | 242/324 | 296/324 | 4811.17 | 4342.37 |
| N54 | **23.46** | 23.83 | 21/22 | 22/22 | 783.00 | 777.00 |

From Table 1 it also can be seen that the heuristic increases the aggregate duration of the tours by increasing the average durations. This behavior is expected because the initial set of tours are optimized to have the shortest durations using the blocking problem. Any changes to those set of tours would cause the overall duration to increase. It is worth noting that this increase is very small. The other key performance indicators summarized in Table 2 (i.e. the aggregate minutes exceeding the 75-minute goal and the ratio of the number of tours exceeding the goal over the number of tours) also show that the heuristic performs well. The heuristic only negatively affects the percentage of the tours exceeding goal KPI in a meaningful way. However, by making modifications to the weight calculation procedure in step 2-2 of the algorithm, we can improve this KPI at the cost of worsening some other. Overall the algorithm does well in balancing as it reduces the standard deviation by up to 47%. It also reduces the minutes exceeding the goal by up to 44%.

## CONCLUSION

In this paper, a Mixed Integer Model (MIP) formulation of the optimal balanced bus scheduling/blocking problem is proposed. Solving this formulation yields a balanced schedule which corresponds to the optimum fleet size. Due to the complexity of the problem, for larger instances, a matching based heuristic algorithm is proposed. The proposed method can obtain leading solutions in significantly shorter time for any size problem. However, using the Xpress MIP solver for solving the proposed MIP is recommended whenever the problem size is tiny so that the solver can find the optimal solution.

To illustrate the performance of the balancing step of the heuristic, it is depicted that how the heuristic effects some statistics and key performance indicators of the tours such as the standard deviation of the tours (SD) and aggregate minutes exceeding the set goal. Statistic measures of the 11 test cases show that the balancing step of the heuristic (the matching phase) always decreases the aggregate minutes exceeding the goal. In fact, it reduces those by up to 44% in some scenarios. It is also indicated that in all scenarios related to HCPSS, the standard deviation of the tours decreases. This reduction in the standard deviation is up to an astonishing 47%.

This study opens the door to research related to balancing of school bus tours. Some possible future research directions could be: performing the balancing using more general methods such as min-cost flow; considering more specific versions of the balanced school bus scheduling problem such as heterogeneous fleet, mixed load, and multiple depots; and incorporation of uncertainties in balanced scheduling.